\documentclass[12pt,fullpage]{article}

\usepackage{amssymb}
\usepackage{amsmath}
\usepackage{epsfig}
\usepackage{caption}
\usepackage{psfrag}
\usepackage{array}

\input epsf

\numberwithin{equation}{section}

\newtheorem{Lemma}{Lemma}[section]
\newtheorem{Cor}[Lemma]{Corollary}
\newtheorem{Th}[Lemma]{Theorem}
\newtheorem{Thm}{Theorem}
\newtheorem{Prop}[Lemma]{Proposition}
\newtheorem{Def}[Lemma]{Definition}
\newtheorem{Rk}[Lemma]{Remark}

\newcommand{\be}{\begin{equation}}
\newcommand{\ee}{\end{equation}}
\newcommand{\baa}{\begin{array}}
\newcommand{\eaa}{\end{array}}
\newcommand{\ba}{\begin{eqnarray}}
\newcommand{\ea}{\end{eqnarray}}
\newcommand{\ds}{\displaystyle}

\def\R{{\mathbb R}}

\def\epsilon{\varepsilon}
\def\trait (#1) (#2) (#3){\vrule width #1pt height #2pt depth #3pt}
\def\fin{\hfill\trait (0.1) (5) (0) \trait (5) (0.1) (0) \kern-5pt \trait (5) (5) (-4.9) \trait (0.1) (5) (0)}

\begin{document}

\date{}
\title{\bf{Robustness for a Liouville type theorem in exterior domains}}
\author{Juliette Bouhours$^{\hbox{\small{1}}}$\\
\\
\footnotesize{$^{\hbox{1}}$UPMC Univ Paris 06, UMR 7598 , Laboratoire Jacques-Louis Lions, F-75005, Paris, France}\\ \footnotesize{e-mail adress: juliette.bouhours@ljll.math.upmc.fr}}

\maketitle

\begin{abstract}
We are interested in the robustness of a Liouville type theorem for a reaction diffusion equation in exterior domains. Indeed H. Berestycki, F. Hamel and H. Matano (2009) proved such a result as soon as the domain satisfies some geometric properties. We investigate here whether their result holds for perturbations of the domain. We prove that as soon as our perturbation is close to the initial domain in the $C^{2,\alpha}$ topology the result remains true while it does not if the perturbation is not smooth enough.\\

\noindent{2010 \em{Mathematics Subject Classification}:  35K57, 35B51,35B53.}\\
\noindent{{\em Keywords:} elliptic equation, Liouville type result, obstacle, maximum principle.}

\end{abstract}


\vskip 20pt
\section{Introduction and main results}\label{intro}
\subsection{Problem and motivations}
This paper investigates the exterior domain problem:
\be\label{problem}
\begin{cases}
-\Delta u=f(u) &\text{in } \R^N\backslash K,\\
\partial_\nu u=0 &\text{on }\partial K,\\
0<u\leq1 &\text{in }  \R^N\backslash K,\\
u(x)\to1\quad \text{as } |x|\to+\infty &\text{uniformly in } x\in\R^N\backslash K,
\end{cases}\ee
where $K$ is a compact set of $\R^N$, $f$ is a bistable non-linearity. 
\vspace{0.1cm}\\
This problem is motivated by the construction of generalized transition fronts for the associated parabolic problem
\be\label{parabolic}\begin{cases}
u_t-\Delta u=f(u) &\text{for all } t\in\R, \quad x\in\R^N\backslash K,\\
\partial_\nu u=0 &\text{for all }t\in\R,\quad x\in\partial K,
\end{cases}\ee
such that 
$$\underset{x\in\R^N\backslash K}{\sup}|u(t,x)-\phi(x_1-ct)|\to0 \text{ as } t\to-\infty,$$
where $\phi$ is a planar traveling wave connecting 1 to 0, i.e
\be\label{TWphi}\begin{cases}
-\phi''-c\phi'=f(\phi) \quad\text{in } \R,\\
\phi(-\infty)=1,\quad \phi(+\infty)=0.
\end{cases}\ee
It is proved in \cite{BHM} that the unique solution of \eqref{parabolic} converges toward a solution of \eqref{problem} as $t\to+\infty$. Thus problem \eqref{problem} determines whether there is a complete invasion or not, that is whether $u(t,x)\to1$ as $t\to+\infty$ \underline{for all} $x\in \R^N\backslash K$. More precisely, complete invasion is shown to hold if and only if \eqref{problem} has no solution different from 1. In \cite{BHM}, Berestycki, Hamel and Matano have shown that if $K$ is star-shaped or directionally convex the unique solution of \eqref{problem} is 1 (see at the end of this section for precise definitions of star-shaped or directionally convex domain). The present paper examines under which conditions this Liouville type theorem is robust under perturbations of the domain. This is shown here to strongly depend on the smoothness of the perturbations that are considered. Indeed our main result is to show that it is true for $C^{2,\alpha}$ perturbations but not for $C^0$ ones. This is stated precisely in the next section. We leave as an open problem to determine what is the optimal space of regularity of the perturbation for which the result remains true.
\vspace{0.2cm}\\
In this paper $f$ is assumed to be a $C^{1,1}([0,1])$ function such that 
\begin{subequations}\label{bistable}
\be
f(0)=f(1)=0, \quad f'(0)<0, \quad f'(1)<0,
\ee
and there exists $\theta\in(0,1)$ such that,
\be
f(s)<0 \quad \forall s\in(0,\theta),\quad f(s)>0\quad \forall s\in(\theta,1).
\ee
\end{subequations}
Moreover we suppose that $f$ satisfies the following positive mass property,
\begin{equation}\label{intf}
\int_0^1f(\tau)d\tau>0.
\end{equation}

\noindent Before stating the main results, let explain what we mean by star-shaped or directionally convex obstacles.
\begin{Def}\label{starshaped}{\rm $K$ is called star-shaped, if either $K=\emptyset$, or there is $x\in\displaystyle{\mathop{K}^{\circ}}$ such that, for all $y\in\partial K$ and $t\in[0,1)$, the point $x+t(y-x)$ lies in $\displaystyle{\mathop{K}^{\circ}}$ and $\nu_K(y)\cdot(y-x)\ge 0$, where $\nu_K(y)$ denotes the outward unit normal to $K$ at $y$.}
\end{Def}
\begin{Def}\label{dirconv}{\rm
K is called directionally convex with respect to a hyperplane $P$ if there exists a hyperplane $P=\{x\in\mathbb{R}^N, x\cdot e=a\}$ where $e$ is a unit vector and $a$ is some real number, such that
\begin{itemize}
\item for every line $\Sigma$ parallel to $e$ the set $K\cap\Sigma$ is either a segment or empty,
\item $K\cap P=\pi(K)$ where $\pi(K)$ is the orthogonal projection of $K$ onto $P$.
\end{itemize}}
\end{Def}


\subsection{ Main results}\label{mainresults}
Our main result is the following Theorem
\begin{Th}\label{constant} 
Let $(K_\epsilon)_{0<\epsilon\leq1}$ be a family of $C^{2,\alpha}$ compact sets of $\R^N$, for some $\alpha>0$. Assume that $K_\epsilon\to K$ for the $C^{2,\alpha}$ topology as $\epsilon\to 0$, and $K$ is either star-shaped or directionally convex with respect to some hyperplane P. Then there exists $\epsilon_0>0$ such that  for all $0<\epsilon<\epsilon_0$, the unique solution of (\ref{problem}) is $u_{\epsilon}\equiv1$
\end{Th}
This theorem means that for obstacles that are compact sets in $\R^N$ and close enough (in the $C^{2,\alpha}$ sense) to some star-shaped or directionally convex domains, the unique solution of (\ref{problem}) is the constant 1. And thus a sufficient condition for the Liouville theorem to be robust under perturbation is the $C^{2,\alpha}$ convergence. On the other hand one can prove that the $C^0$ convergence of the perturbation is not enough for the result to stay true. This is stated in the Theorem below.
\begin{Th}\label{counterex}
There exists $(K_\epsilon)_\epsilon$ a family of compact manifolds of $\R^N$ such that $K_\epsilon\to B_{R_0}$ for the $C^0$ topology as $\epsilon\to 0$, and for all $\epsilon>0$ there exists a solution $u_\epsilon$ of \eqref{problem} such that $0<u_\epsilon<1$ in $\R^N\backslash K_\epsilon$.
\end{Th}

\noindent \textbf{Notations}\\
We denote by $B_{R_0}$ the ball of radius $R_0$ centred at 0 in $\R^N$, i.e
$$B_{R_0}:=\left\{x\in\R^N,\: |x|<R_0\right\},$$
and by $B_{r}(x_0)$ the ball of radius $r$ centred at $x_0$ in $\R^N$, i.e
$$B_r(x_0):=\left\{x\in\R^N,\: |x-x_0|<r\right\},$$

\begin{Rk}[$C^0$ or $C^{2,\alpha}$ convergence]\label{C2conv}
When we write $K_\varepsilon\to K$ for the $X$ topology we mean that for each $x_0$ in $\partial K$, and for some $r>0$ such that $\partial K_\epsilon\cap B_r(x_0)\neq\emptyset$ there exists a couple of functions $\psi_\epsilon$ and $\psi$ defined on $B_r(x_0)$, parametrization  of $K_\epsilon$ and $K$, such that, $\psi_\epsilon\in X(B_r(x_0))$ and $\psi\in X(B_r(x_0))$ with $\lVert\psi_\epsilon-\psi\rVert_{X(B_r(x_0))}\to0$ as $\epsilon\to0$.\\
For more details about the $C^{2,\alpha}$ topology one can look at \cite{GT}, chapter 6.
\end{Rk}

\noindent Before proving the previous statements, let give some examples of domains $(K_\epsilon)_\epsilon$ and $K$ to illustrate our results.


\subsection{Examples of domains}
We assume that $N=2$ and we construct two families of obstacles; one which converges to a star-shaped domain and the other which converges to a directionally convex domain. The black plain line represents the limit $K$ and the thin parts represent the small perturbations (of order $\epsilon$).
\vspace{0.2cm}\\
\begin{minipage}[b]{0.45\linewidth}
\centering
\includegraphics[scale=0.3]{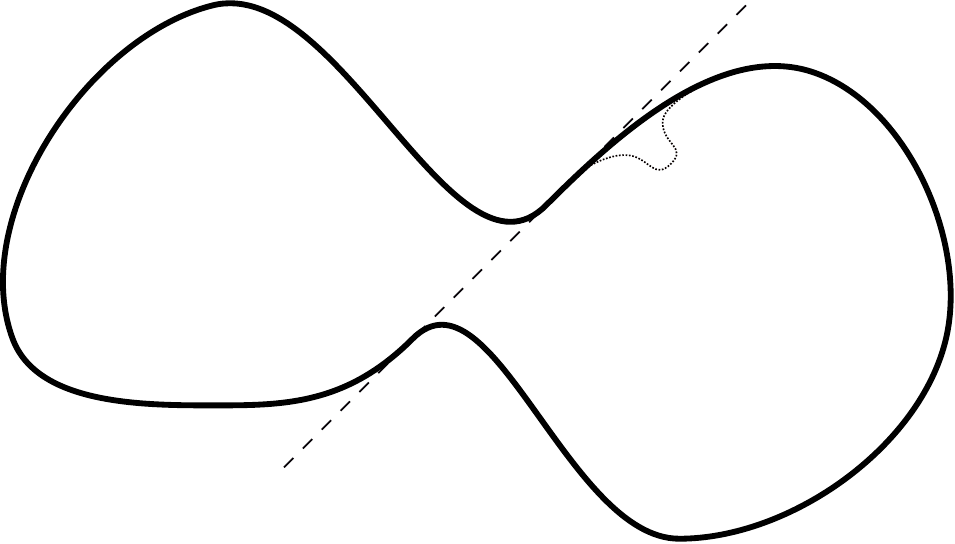}
\captionof{figure}{\footnotesize{Obstacle that converges toward a star-shaped domain}} 
\label{starshapedC2}
\end{minipage}
\hspace{0.5cm}
\begin{minipage}[b]{0.45\linewidth}
\centering
\includegraphics[scale=0.2]{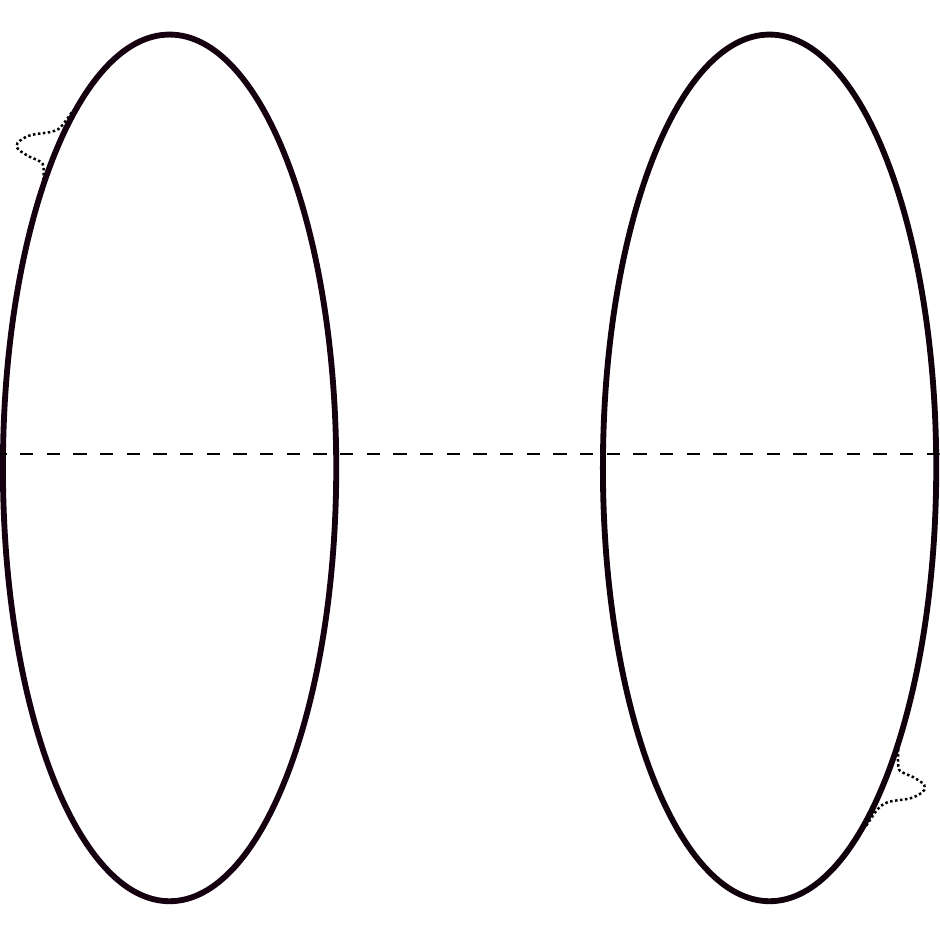}
\captionof{figure}{\footnotesize{Obstacle that converges toward a directionally convex domain}}
\label{dirconvC2}
\end{minipage}
\vspace{0.2cm}\\
The long-dashed  lines are used during the construction of $K$. For the star-shaped domain, it is on this line that we could find the center(s) of the domain (i.e the point $x$ in Definition \ref{starshaped}). For the directionally convex domain it represents the hyperplane $P$.  We can clearly see that for all $\epsilon>0$, $K_\epsilon$ is not star-shaped for figure \ref{starshapedC2} and not  directionally convex for figure \ref{dirconvC2}. One needs to be careful on the shape of the perturbations. Indeed in figure \ref{ellipseC0} below $K_\epsilon$ converges to an ellipse as $\epsilon\to0$ but the convergence of $K_\epsilon$ is not $C^{2,\alpha}$ (see section \ref{Discussion} for more details) but only $C^0$ which is not enough for the Liouville theorem to remain valid.
\vspace{0.2cm}
\begin{center}
\includegraphics[scale=0.7]{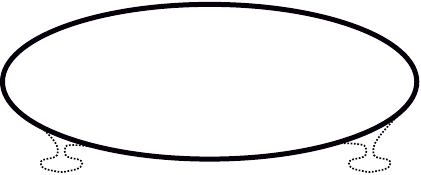}
\captionof{figure}{Obstacles converging only in the $C^0$ topology} 
\label{ellipseC0}
\end{center}
\vspace{0.2cm}
We will prove Theorem \ref{constant} in section \ref{proof} below and Theorem \ref{counterex} in section \ref{Discussion}.

\section{Robustness of the result for $C^{2,\alpha}$ perturbations}\label{proof}
In this section we prove the robustness of the Liouville result when the perturbation is close to a star-shaped or directionally convex domain in the $C^{2,\alpha}$ topology. To prove Theorem \ref{constant}, we will need the following Proposition:
\begin{Prop}\label{unifconv}
For all $0<\delta<1$, there exists $\ds{R=R_\delta}$ such that, if $u_\epsilon$ is a solution of \eqref{problem} with $K_\epsilon$, then $u_\varepsilon(x)\geq 1-\delta$ for all $\lvert x\rvert\geq R$ and for all $0<\varepsilon<1$.
\end{Prop}
This proposition means that $u_\epsilon$ converges toward 1 as $\lvert x\rvert\to+\infty$ uniformly in $\epsilon$.\\
Let first admit this result and prove Theorem \ref{constant}.


\subsection{Proof of Theorem \ref{constant}}
As $u_\epsilon$ is uniformly bounded for all $\epsilon >0$, using Schauder estimates, we know that up to a subsequence $u_{\epsilon_n}\to u_0$ in $C^2_\textrm{loc}$ as $n\to +\infty$, in the sense that for all $r>0$, 
$$\lVert u_{\epsilon_n}-u_0\rVert_{C^2(B_r\backslash (K_{\epsilon_n}\cup K))}\to0 \text{ as } n\to+\infty$$ and $u_0$ satisfies:
\begin{equation}
\begin{cases} \Delta u_0+f(u_0)=0 & \textrm{in } \mathbb{R}^N\backslash K,\\
		\nu \cdot \nabla u_0=0 &\textrm{on } \partial K.\\
	
\end{cases}
\end{equation}
Using Proposition \ref{unifconv} we get $\underset{\lvert x\rvert \to +\infty} {\lim}u_0(x)=1$. And $K$ is either star-shaped or directionally convex. We now recall the following results from \cite{BHM}:
\begin{Thm}[Theorem 6.1 and 6.4 in \cite{BHM}] \label{liouville}
Let $f$ be a Lipschitz-continuous function in $[0,1]$ such that $f(0)=f(1)=0$ and $f$ is nonincreasing in $[1-\delta,1]$ for some $\delta>0$. Assume that 
\begin{equation}
\forall \hspace{0.1cm}0\leq s<1,\hspace{0.2cm} \int_s^1f(\tau)d\tau>0.
\end{equation}
Let $\Omega$ be a smooth, open, connected subset of $\R^N$ (with $N\geq2$) with outward unit normal $\nu$, and assume that $K=\R^N\backslash\Omega$ is compact. Let $0\leq u\leq1$ be a classical solution of
\begin{equation}\label{laplNeu}
\begin{cases}-\Delta u=f(u) &\textrm{in } \Omega,\\
\nu\cdot\nabla u=0 &\textrm{on } \partial\Omega,\\
u(x)\to1 &\textrm{as } \lvert x \rvert\to+\infty.
\end{cases}
\end{equation}
If $K$ is star-shaped or directionally convex, then 
\begin{equation}
u\equiv 1 \textrm{ in } \overline{\Omega}.
\end{equation}
\end{Thm} 
It thus follows that $u_0\equiv1$.
It also proves that the limit $u_0$ is unique and thus $u_\epsilon\to u_0$ as $\epsilon \to 0$ in $C^2_\text{loc}$ (and not only along a subsequence). 
\vspace{0.2 cm}\\
Now we need to prove that there exists $\epsilon_0>0$ such that $u_\epsilon\equiv 1$ for all $0<\epsilon<\epsilon_0$.
Let assume that for all $\epsilon>0$, $u_\varepsilon\not\equiv1$ in $\Omega_\epsilon=\R^N\backslash K_\epsilon$. Then there exists $\ds{x_0 \in \overline{\Omega}_\varepsilon}$ such that $\ds{u_\varepsilon(x_0)=\underset{x\in\overline{\Omega_\varepsilon}}{\min } \hspace{0.1cm} u_\varepsilon(x)}<1$ and $x_0$, which depends on $\epsilon$, is uniformly bounded with respect to $\epsilon$ (using Proposition \ref{unifconv}). Assume that  $u_\epsilon(x_0)>\theta$, as $u_\epsilon$ is a solution of (\ref{problem}), the Hopf lemma yields that,
\begin{equation*}
\textrm{if } x_0\in\partial K_\varepsilon \textrm{ then }  \frac{\delta u_\epsilon}{\delta \nu}(x_0)<0,
\end{equation*}
which is impossible due to Neumann boundary conditions. Hence $x_0 \in \Omega_\varepsilon$ and 
\begin{equation*}
-\Delta u_\varepsilon(x_0)=f(u_\varepsilon(x_0))>0,
\end{equation*}
which is impossible since $x_0$ is a minimizer. So, for all $0<\epsilon<1$,
\begin{equation*}
0\leq \underset{x\in\overline{\Omega_\varepsilon}}{\min } \hspace{0.1cm} u_\varepsilon(x)\leq \theta,
\end{equation*}
But $\underset{x\in\overline{\Omega_\epsilon}}{\min } u_\epsilon\to1$ as $\epsilon\to0$ by Proposition \ref{unifconv} and the local uniform convergence to $u_0\equiv1$, which is a contradiction.
Thus there exists  $\varepsilon_0$ such that for all $\varepsilon<\varepsilon_0$, $u_\varepsilon \equiv 1$. \fin


\subsection{Proof of Proposition \ref{unifconv}}
We will now prove Proposition \ref{unifconv}, using the following lemma:
\begin{Lemma}\label{omega}
There exists $\omega=\omega(r)$ with $r\in \mathbb{R}^+$ such that 
\begin{equation}\label{lemmaomega}
\begin{cases}
-\omega''(r)=f(\omega(r)), & \forall r\in \R^+_*,\\
\omega(0)=0, \hspace{0.2cm} \omega'(0)>0,\\
 \omega'>0,\hspace{0.2cm} 0<\omega<1 & \textrm{in } \mathbb{R}_+^*,\\
 \underset{r \to+\infty}{\lim}\omega(r)=1.
 \end{cases}
 \end{equation}
\end{Lemma}
This is a well known result. In deed, by a shooting argument, if $\omega$ is a solution of the initial value problem
\begin{equation*}\label{edo}
\begin{cases} -\omega''=f(\omega) & \textrm{in } (0,+\infty), \\
		\omega(0)=0,\\
		\omega'(0)=\sqrt{2F(1)},
	
\end{cases}
\end{equation*}
where $F(1)=\int_0^1f(s)ds$, it is easily seen that $\omega$ is also a solution of \eqref{lemmaomega}.
\vspace{0.3 cm}\\
\noindent{\bf{Proof of Proposition \ref{unifconv}.}}
Now we introduce a function $f_\delta$ defined in $[0,1-\frac{\delta}{2}]$, satisfying the same bistability hypothesis as $f$ but such that 
\begin{itemize}
\item[$\cdot$] $f_\delta\leq f$ in $[0,1-\frac{\delta}{2}]$, 
\item[$\cdot$] $f_\delta=f$ in $[0,1-\delta]$,  
\item[$\cdot$] $f_\delta(1-\frac{\delta}{2})=0$.
\end{itemize}
Notice that $\int_0^{1-\frac{\delta}{2}}f_\delta(z)dz>0$ for $\delta$ small.
Using the same arguments than in Lemma \ref{omega} there exists $\omega=\omega_\delta$ such that 
\begin{equation}\label{TWdelta}
\begin{cases} -\omega_\delta''(x)=f_\delta(\omega_\delta(x)) & \textrm{in } (0,+\infty),\\
		\omega_\delta(0)=0,\hspace{0.2cm}\omega_\delta(+\infty)=1-\frac{\delta}{2},\\
		0<\omega_\delta<1-\frac{\delta}{2}& \textrm{in } (0,+\infty) ,\\
		\omega_\delta'>0 & \textrm{in } (0,+\infty).
	
\end{cases}
\end{equation}
As $K_\epsilon$ is a compact set of $\R^N$ converging to a fix compact set $K$, there exists $R_0$ such that $K_\epsilon\subset B_{R_0}$ for all $\epsilon>0$.\\
Next, for any $R>R_0$ let consider $z(x)=\omega_\delta (\lvert x\rvert-R)$, for every $\lvert x\rvert\geq R$. One gets:
\begin{equation}\label{subsol}
-\Delta z<f(z) \text{ in } \R^N\backslash B_R.
\ee 
We want to prove that
$$\omega_\delta(\lvert x \rvert -R_0)< u_\varepsilon(x), \hspace{1cm} \forall x\in \mathbb{R}^N, \lvert x \rvert\geq R_0.$$
We know from (\ref{problem}) that $u_\varepsilon(x)\to 1$ as $\lvert x\rvert\to+\infty$. Hence there exists $A=A(\varepsilon)>0$ such that $\ds{u_\varepsilon(x)\geq1-\frac{\delta}{3}>\omega_\delta(|x|-A)}$, for all $\lvert x\rvert\geq A$.
Consider 
\be\label{Rbar}
\overline{R}=\inf \left\{R\geq R_0; u_\varepsilon(x)> \omega_\delta(\lvert x\rvert-R)\textrm{, for all } \lvert x\rvert\geq R\right\}.
\ee
As $\overline{R}\geq R_0$ and $K_\epsilon\subset B_{R_0}$, $u_\epsilon$ is always defined in $\ds{\left\{|x|>\overline{R}\right\}}$. One will prove that $\overline{R}=R_0$. As $\omega_\delta$ is increasing, we know that $$\forall R\geq A \hspace{0.3 cm}u_\varepsilon(x)\geq \omega_\delta(\lvert x\rvert -R),\hspace{0.3 cm} \forall \lvert x\rvert\geq R.$$
Hence $\overline{R}\leq A$.\\
Assume that $\overline{R}>R_0$. Then there are two cases to study:
\begin{itemize}
\item either $\underset{\lvert x\rvert>\overline{R}}{\inf}\Big\{u_\varepsilon(x)-\omega_\delta (\lvert x\rvert-\overline{R})\Big\}>0$, (1)
\item or $\underset{\lvert x\rvert>\overline{R}}{\inf}\Big\{u_\varepsilon(x)-\omega_\delta (\lvert x\rvert-\overline{R})\Big\}=0$. (2)
\end{itemize}
In the first case (1), one gets $u_\varepsilon(x)>\omega_\delta (\lvert x\rvert-\overline{R})$ for all $ \lvert x\rvert>\overline{R}$. As $\nabla u_\varepsilon$ and $\omega_\delta'$ are bounded, there exists $R^*<\overline{R}$ such that $u_\varepsilon(x)\geq\omega_\delta (\lvert x\rvert-R^*)$ for all $|x|>R^*$. This contradicts the optimality of $\overline{R}$.
\vspace{0.4cm}\\
In the second case (2), there necessarily exists $x_0$ with $\lvert x_0\rvert>\overline{R}$ such that $u_\varepsilon (x_0)=\omega_\delta(\lvert x_0\rvert-\overline{R})$. Let $v(x)=u_\epsilon(x)-\omega_\delta(|x|-\overline{R})$, for all $|x|>\overline{R}$. As $u_\epsilon$ is a solution of (\ref{problem}) and using (\ref{subsol}), $v$ satisfies:
\be\label{v}\begin{cases}
-\Delta v>c(x)v &\textrm{in }\left\{|x|>\overline{R}\right\},\\
v>0 &\textrm{on } \left\{|x|=\overline{R}\right\},
\end{cases}
\ee 
where $c$ is a bounded function. From the definition of $\overline{R}$, $v(x)\geq 0$, for all $\lvert x \rvert \geq \overline{R}$. But there exists $x_0$ such that $\lvert x_0\rvert>\overline{R}$ and $v(x_0)=0$ which implies that $v(\cdot)\equiv 0$. This is impossible because $v(\cdot)>0$, for all $\lvert x \rvert=\overline{R}$.
\vspace{0.2cm}\\
Then $\overline{R}=R_0$ which does not depend on $\varepsilon$ and
$$\forall  \lvert x\rvert \geq R_0 \hspace{0.3cm} u_\varepsilon(x)\geq\omega_\delta(\lvert x\rvert-R_0).$$
As $\omega_\delta(x)\to 1-\frac{\delta}{2}$ as $\lvert x\rvert\to+\infty$, there exists $\hat{R}$, independent of $\epsilon$, such that for all $\lvert x \rvert>\hat{R}+R_0$, $u_\epsilon(x)>\omega_\delta(\lvert x\rvert-R_0)\geq 1-\delta$. 
One has proved Proposition \ref{unifconv}. \fin


\section{Counter example in the case of $C^0$ perturbations}\label{Discussion}
Until now we have assumed that  $K_\varepsilon\to K$ in $C^{2,\alpha}$, in order to use the Schauder estimates and ensure the convergence of $u_\epsilon$ as $\epsilon\to 0$ . One can wonder if we can weaken this hypothesis, i.e would the $C^0$ or $C^1$ convergence be enough?\\
We prove that $C^0$ perturbations are not smooth enough for the Liouville result to remain true. 


\subsection{Construction of a particular family of $C^0$ perturbations}
In this subsection we construct a family of obstacles that are neither star-shaped nor directionally convex but converges uniformly to $B_{R_0}$ which is convex. We want to prove that for all $\varepsilon \in ]0,1]$ there exists a solution of \eqref{problem} which is not identically equal to 1. To do so we will use the counterexample of section 6.3 in \cite{BHM}.
\begin{center}
\includegraphics[scale= 0.8]{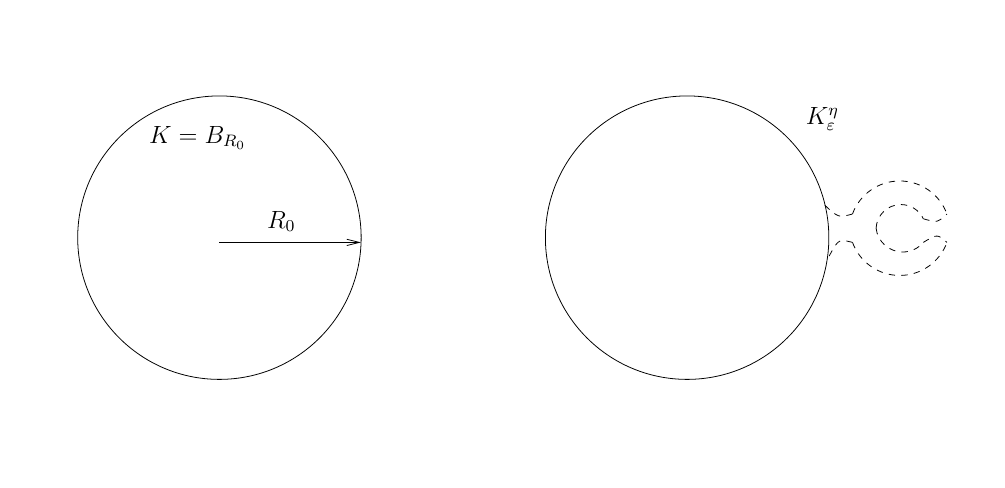}
\captionof{figure}{Liouville counterexample}\label{obstacle} 
\label{liouv}
\end{center}
Zooming on the dashed part:
\begin{center}
\includegraphics[scale=0.9]{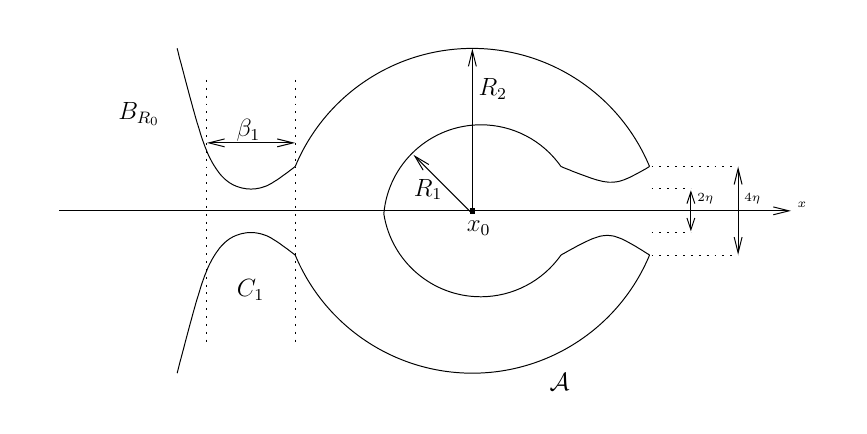}
\captionof{figure}{Zoom on the perturbation}\label{zoom}
\label{obsC0zoom}
\end{center}
We consider an obstacle $K_1=K_1^\eta$ (see figure \ref{obstacle} and \ref{zoom}), such that:
\begin{equation}\label{defobsC0}
\begin{cases}
\big(\mathcal{A}\cap\{x;x_1\leq x_1^0\}\big)\cup B_{R_0}\cup C_1 \subset K_1^\eta,\\
\mathcal{A}\cap\big\{x;x_1>x_1^0, \lvert x'\rvert>2\eta\big\} \subset K_1^\eta,\\
K_1^\eta\subset \Big(\mathcal{A}\cap\big\{x;x_1> x_1^0,\lvert x'\rvert>\eta\big\}\Big)\cup B_{R_0}\cup\Big(\mathcal{A}\cap\big\{x;x_1\leq x_1^0\big\}\Big)\cup C_1.
\end{cases}
\end{equation}
where $x'=(x_2,...,x_N)$ and $\mathcal{A}=\{x: R_1\leq \lvert x-x^0\rvert\leq R_2\}$, $R_0$, $R_1<R_2$, are three positive constants, $x^0=(x^0_1,0,0,...,0)$ is the center of the annular region $\mathcal{A}$ with $x_1^0=R_0+R_2+\beta_1$, $C_1$ is some corridor that links smoothly $\mathcal{A}$ and $B_{R_0}$ which length is $\beta_1$ and $\eta>0$, small enough.\\
The family $(K_\epsilon)$ is constructed by downsizing $K_1$ such that for all $0<\epsilon<1$, $\mathcal{A}_\epsilon$ stays an annular region, with 
\begin{itemize}
\item[$\cdot$] $x_0^\epsilon=(R_0+R_2^\epsilon+\beta_\epsilon, 0)\in\R^N$ converging to $(R_0, 0)\in\R^N$, 
\item[$\cdot$] $R_1^\epsilon=\epsilon R_1$, $R_2^\epsilon=\epsilon R_2$,
\item[$\cdot$]  $\beta_\epsilon$ converging to 0 as $\epsilon\to0$. 
\end{itemize}
We have the following lemma.
\begin{Lemma}\label{C0notC1}
$K_\epsilon\to K$ for the $C^0$ topology as $\epsilon\to 0$ but not for the $C^1$ topology.
\end{Lemma}
This Lemma is easily proved using smooth parametrization of $B_{R_0}$ and $K_1$ and noticing that for all $\epsilon>0$ there exists a point on the boundary of the perturbation that has an outward unit normal orthogonal to $e_1=(1,0,...,0)$.
%


\subsection{Existence of a non constant solution $u_\epsilon$ of (\ref{problem})}
We want to prove that for all $0<\epsilon<1$ there exists a solution $0<u_\epsilon<1$ of 
\begin{equation}\label{probeps}
\begin{cases}
-\Delta u_\epsilon=f(u_\epsilon) & \textrm{in } \R^N\backslash K_\epsilon^\eta=\Omega_\epsilon^\eta,\\
 \nu\cdot\nabla u_\epsilon=0 &\textrm{on } \partial K_\epsilon^\eta=\partial\Omega_\epsilon^\eta,\\
 u_\epsilon(x)\to 1 \textrm{ as } \lvert x\rvert\to+\infty.
 \end{cases}
 \end{equation}
We will follow the same steps as in \cite{BHM}, section 6. First, let notice that it is enough to find $\omega\not\equiv1$ solution of 
\begin{equation}\label{probBR}
\begin{cases}
-\Delta \omega=f(\omega) & \textrm{in } B_R\backslash K_\epsilon^\eta,\\
\nu\cdot\nabla\omega=0 &\textrm{on } \partial K_\epsilon^\eta,\\
\omega=1 &\textrm{on } \partial B_R,
\end{cases}
\end{equation}
for some $R>0$ large enough such that $K_\epsilon^\eta\subset B_R$.\\
Indeed then $\omega$ extended by 1 outside $B_R$ is a supersolution of (\ref{probeps}) and one can define:
\begin{equation}\label{psi}
\psi(x)= \begin{cases} 	0 &\textrm{ if } \{\lvert x\rvert<R\}\backslash K_\epsilon^\eta,\\
				U(\lvert x \rvert-R) &\textrm{ if } \lvert x \rvert \geq R,
\end{cases}
\end{equation}
where $U:\R^+\to(0,1)$ satisfies $U''+f(U)=0$ in $\R^*_+$, $U(0)=0$, $U'(\xi)>0$ $\forall$ $\xi\geq0$, $U(+\infty)=1$. It exists as soon as (\ref{intf}) is satisfied (see Lemma \ref{omega}). As $U(\lvert\cdot\rvert-R)$ is a subsolution, $\psi$ is a subsolution of (\ref{probeps}).\\
Hence there exists a solution $\psi<u_\epsilon<\omega$ of (\ref{probeps}). If we prove that $\omega\not\equiv 1$ then $0<u_\epsilon<1$ (with the maximum principle).
\vspace{0.2cm}\\
Now consider our problem (\ref{probBR}) and replace $\omega$ by $v=1-\omega$. The problem becomes
\begin{equation}\label{probBRv}
\begin{cases}
-\Delta v=-f(1-v)=g(v) & \textrm{in } B_R\backslash K_\epsilon^\eta,\\
\nu\cdot\nabla v=0 &\textrm{on } \partial K_\epsilon^\eta,\\
v=0 &\textrm{on } \partial B_R.
\end{cases}
\end{equation}
Using exactly the same arguments as in \cite{BHM} one proves that, if we consider:
\begin{equation}\label{condobsC0}
v_0(x)=
\begin{cases}1 & \hspace{0.53cm} \textrm{if } x\in B_{R_2}(x^0)\backslash K_\epsilon^\eta\cap\Big\{x; x_1-x_1^0\leq \frac{2R_1+R_2}{3}\Big\},\\
\frac{3}{R_2-R_1}(\frac{R_1+2R_2}{3}-(x_1-x_1^0)) & \begin{array}{clcr}
 &\textrm{if } x\in B_{R_2}(x^0)\backslash K_\epsilon^\eta\\
 &\cap\Big\{x;\frac{2R_1+R_2}{3}\leq x_1-x_1^0\leq \frac{R_1+2R_2}{3}\Big\},
 \end{array}\\
0 &\begin{array}{clcr} 
&\textrm{if } x\in \Big[B_R\backslash\big(B_{R_2}(x^0)\cup C_\epsilon\cup B_{R_0}(0)\big)\Big]\\
&\cup \Big[B_{R_2}(x^0)\backslash K_\epsilon^\eta \cap\big\{x, x_1-x_1^0\geq \frac{R_1+2R_2}{3}\big\}\Big],
\end{array}
\end{cases}
\end{equation}
then for $\eta>0$ small enough, there exists $v\in H^1(B_R\backslash K_\epsilon^\eta)\cap\{v=0 \textrm{ on } \partial B_R\}=\overline{H}^1_0$, $\delta>0$ such that $\lVert v-v_0\rVert_{H^1}<\delta$ and $v$ is a local minimizer of the associated energy functional in $\overline{H}^1_0$. For more clarity we will give the main steps of the proof but for details and proofs see \cite{BHM}, section 6.3. 
\vspace{0.2cm}\\
We introduce the energy functional in a domain $D$:
\begin{equation}\label{energy}
J_D(\omega)=\int_D\big\{\frac{1}{2}\lvert\nabla\omega\rvert^2-G(\omega)\big\}dx,
\end{equation}
defined for functions of $H^1(D)$, where
\begin{equation}\label{funcG}
G(t)=\int_0^tg(s)ds,
\end{equation}
$g$ defined in \eqref{probBRv}. Using Proposition 6.6 in \cite{BHM} one gets the following Corollary
\begin{Cor}\label{energyD0}
In $B_{R_1}(x^0)$, $v_0\equiv 1$ is a strict local minimum of $J_{B_{R_1}(x^0)}$ in the space $H^1(B_{R_1}(x^0))$. More precisely, there exist $\alpha>0$ and $\delta>0$ for which
\begin{equation}\label{minD0}
J_{B_{R_1}(x^0)}(v)\geq J_{B_{R_1}(x^0)}(v_0)+\alpha \lVert v-v_0\rVert^2_{H^1(B_{R_1}(x^0))},
\end{equation}
for all $v\in H^1(B_{R_1}(x^0))$ such that $\lVert v-v_0\rVert^2_{H^1(B_{R_1}(x^0))}\leq\delta$.
\end{Cor}
And then using Proposition 6.8 of \cite{BHM} and Corollary \ref{energyD0} one gets the following Corollary
\begin{Cor}\label{energymin}
There exist $\gamma>0$ and $\eta_0>0$ (which depend on $\epsilon$) such that for all $0<\eta<\eta_0$ and $v\in \overline{H}^1_0$ such that $\lVert v-v_0\rVert^2_{B_R\backslash K_\epsilon^\eta}=\delta$, then
\begin{equation*}
J_{B_R\backslash K_\epsilon^\eta}(v_0)<J_{B_R\backslash K_\epsilon^\eta}(v)-\gamma.
\end{equation*}
\end{Cor}
The proof of this corollary relies on the existence of a channel of width  of order $\eta>0$ opening on the interior of the annular region $\mathcal{A}$ (third assumption in \eqref{condobsC0}). This condition cannot be satisfied if the convergence of the obstacle is $C^1$ (see  Lemma \ref{C0notC1}).\\
The functional $J_{B_R\backslash K_\epsilon^\eta}$ admits a local minimum in the ball of radius $\delta$ around $v_0$ in $H^1(B_R\backslash K_\epsilon^\eta)\cap\{v=0 \textrm{ on } \partial B_R\}$. This yields a (stable) solution $v$ of (\ref{probBRv}) for small enough $\eta>0$. Furthermore, provided that $\delta$ is chosen small enough, this solution does not coincide either with 1 or with 0 in $B_R\backslash K_\epsilon^\eta$.\\
We have proved that for all $\varepsilon \in ]0,1]$, $0<u_\varepsilon< 1$.
\vspace{0.5cm}\\
One has proved that $C^0$ convergence of the domain is not sufficient and thus Theorem \ref{counterex}. 
\vspace{0.2cm}\\
\noindent One can conclude that if the perturbation is smooth in the $C^{2,\alpha}$ topology, we still have a Liouville type result for reaction diffusion equation in exterior domain. Whereas one can construct a counterexample of this Liouville result for $C^0$ perturbations. One question that is still open is thus the optimal space of regularity of the perturbation for the results to remain true under perturbation. For instance is the $C^1$ convergence of the perturbation enough to get the result?

\section*{Acknowledgments}
The author thanks Henri Berestycki, Gr\'egoire Nadin and Fran\c{c}ois Hamel for their fruitful discussions and encouragements, Sylvain Arguill\`ere for his geometric insights. This study was supported by the French "Agence Nationale de la Recherche" through the project PREFERED (ANR 08-BLAN-0313) and by the French "R\'egion Ile de France" through the fellowship "Bourse Hors DIM".

\vskip 20pt
\bibliography{bibliothese}

\def\cprime{$'$}
\begin{thebibliography}{1}

\bibitem{BHM}
Henri Berestycki, Fran{\c{c}}ois Hamel, and Hiroshi Matano.
\newblock Bistable traveling waves around an obstacle.
\newblock {\em Comm. Pure Appl. Math.}, 62(6):729--788, 2009.

\bibitem{GT}
David Gilbarg and Neil~S. Trudinger.
\newblock {\em Elliptic partial differential equations of second order}.
\newblock Classics in Mathematics. Springer-Verlag, Berlin, 2001.
\newblock Reprint of the 1998 edition.

\end{thebibliography}
\bibliographystyle{plain}

\end{document}